\documentclass[12pt]{amsart}
\usepackage{amssymb, mathrsfs}

\long\def\forget#1\forgotten{}
\newcommand{\Impl}{\Rightarrow}
\newcommand{\bi}{\begin{itemize}}
\newcommand{\itm}{\item}
\newcommand{\ei}{\end{itemize}}
\newcommand{\be}{\begin{enumerate}}
\newcommand{\ee}{\end{enumerate}}
\newcommand{\Tau}{\mathrm{T}}
\newcommand{\nin}{\notin}
\newcommand{\cO}{\mathcal{O}}
\newcommand{\Union}{\bigcup}
\newcommand{\sbst}{\subseteq}

\newcommand{\N}{\mathbb{N}}
\newcommand{\seq}[1]{\{#1\}_{n\in\N}}

\newcommand{\scrA}{\mathscr{A}}
\newcommand{\scrB}{\mathscr{B}}

\newcommand{\sm}{\setminus}
\newcommand{\as}{\subseteq^*}\newcommand{\roth}{{[\N]^{\aleph_0}}}
\theoremstyle{plain}
\newtheorem{thm}{Theorem}
\newtheorem{prop}[thm]{Proposition}
\newtheorem{lem}[thm]{Lemma}
\theoremstyle{plain}
\newtheorem{cor}[thm]{Corollary}
\theoremstyle{definition}
\newtheorem{defn}[thm]{Definition}
\theoremstyle{remark}
\newtheorem{rem}[thm]{Remark}
\newcommand{\grp}[1]{\gimel(#1)}
\newcommand{\grpa}{\grp{\scrA}}
\newcommand{\grpb}{\grp{\scrB}}
\newcommand{\sone}{\mathsf{S}_{1}}
\newcommand{\sfin}{\mathsf{S}_\mathrm{fin}}
\newcommand{\ufin}{\mathsf{U}_\mathrm{fin}}

\newcommand{\gone}{\mathsf{G}_{1}}

\newcommand{\gfin}{\mathsf{G}_\mathrm{fin}}

\newcommand{\cU}{\mathcal{U}}
\newcommand{\cV}{\mathcal{V}}
\newcommand{\cW}{\mathcal{W}}
\newcommand{\cF}{\mathcal{F}}
\newcommand{\cP}{\mathcal{P}}

\title[Partition relations for Hurewicz selections]{Partition relations for Hurewicz-type selection hypotheses}

\author{Nadav Samet}
\address[Nadav Samet]{Department of Mathematics,
Weizmann Institute of Science, Rehovot 76100, Israel
}
\curraddr{Google Ireland Ltd., Gordon House, Barrow Street, Dublin 4, Ireland}
\email{thesamet@gmail.com}

\author{Marion Scheepers}
\address[Marion Scheepers]{Department of Mathematics, Boise State University, Boise, ID 83725, USA}
\email{marion@diamond.boisestate.edu}

\author{Boaz Tsaban}
\address[Boaz Tsaban]{Department of Mathematics, Bar-Ilan University, Ramat-Gan 52900, Israel;
and
Department of Mathematics, Weizmann Institute of Science, Rehovot 76100, Israel
}
\email{tsaban@math.biu.ac.il}

\begin{document}
\begin{abstract}
We give a general method to reduce Hurewicz-type selection hypotheses into
standard ones. The method covers the known results of this kind
and gives some new ones.

Building on that, we show how to derive Ramsey theoretic
characterizations for these selection hypotheses.
\end{abstract}

\subjclass[2000]{
05C55, 
05D10, 
54D20 
}

\keywords{Ramsey theory of open covers, selection principles, Hurewicz covering property, $\tau$-covers}

\maketitle

\section{Introduction}
In \cite{Menger24}, Menger introduced a hypothesis which
generalizes $\sigma$-compact\-ness of topological spaces. Hurewicz
\cite{Hure27} proved that Menger's property is equivalent to a
property of the following type.
\begin{description}
\item[$\sfin(\scrA,\scrB)$]
For each sequence $\seq{\cU_n}$ of members of $\scrA$,
there exist finite subsets $\cF_n\sbst\cU_n$, $n\in\N$, such that $\Union_n\cF_n\in\scrB$.
\end{description}
Indeed, Hurewicz observed that $X$ has Menger's property if, and only if,
$X$ satisfies $\sfin(\cO,\cO)$, where $\cO$ is the collection of
all open covers of $X$.
Motivated by a conjecture of Menger, Hurewicz \cite{Hure27} introduced
a hypothesis of the following type.
\begin{description}
\item[{$\ufin(\scrA,\scrB)$}] For each sequence
$\{\cU_{n}\}_{n\in\mathbb{N}}$ of elements of $\scrA$ which
do not contain a finite subcover, there exist finite (possibly
empty) subsets $\mathcal{F}_{n}\subseteq\cU_{n}$,
$n\in\mathbb{N}$, such that
$\{\bigcup\mathcal{F}_{n}:n\in\mathbb{N}\}\in\scrB$.
\end{description}
Hurewicz was interested in $\ufin(\cO,\Gamma)$, where $\Gamma$ is
the collection of all open $\gamma$-covers of $X$. ($\cU$ is a \emph{$\gamma$-cover} of
$X$ if it is infinite, and each $x\in X$ is an element in all but finitely many
members of $\cU$.)

While the Hurewicz-type selection hypotheses $\ufin(\scrA,\scrB)$
are standard notions in the field of selection principles,
they are less standard in the more general field of infinitary
combinatorics. The reason for that is that the finite subsets are
``glued'' before considering the resulting object. Moreover, the
definition of $\ufin(\scrA,\scrB)$ is somewhat less elegant than
that of $\sfin(\scrA,\scrB)$, and consequently is less convenient
to work with.

$\cU$ is an \emph{$\omega$-cover} of $X$ if $X\nin \cU$, but for each finite $F\sbst X$,
there is $U\in\cU$ such that $F\sbst U$.
One of the main results of \cite{coc7} is the result that
$\ufin(\cO,\Gamma)$ is equivalent to $\sfin(\Omega,\scrB)$ for an
appropriate modification $\scrB$ of $\Gamma$.
A similar result was established in \cite{coc8} for $\ufin(\cO,\Omega)$.
In these papers, these reductions were used to obtain Ramsey theoretic
characterizations of $\ufin(\cO,\Gamma)$ and of $\ufin(\cO,\Omega)$,
respectively.

We generalize these results. As applications,
we reproduce the main results of \cite{coc7},
strengthen the main results of \cite{coc8}, and
obtain standard and Ramsey theoretic equivalents
for a property introduced in \cite{tautau}.

As each cover of a type considered in our Ramsey theoretic results
is infinite, each of our Ramsey theoretic results implies Ramsey's classical
theorem and can therefore be viewed as a structural extension of
Ramsey's Theorem.

\section{Reduction of selection hypotheses}

\subsubsection*{Convention}
To simplify the presentation, by \emph{cover} of $X$ we always mean a \emph{countable} collection $\cU$ of
open subsets of $X$, such that $\bigcup\cU=X$ and $X\notin\cU$. Also,
$\scrA$ and $\scrB$ always denote families of (such) covers of the underlying space $X$.

\begin{defn}
$\scrA$ is \emph{Ramseyan} if for each $\cU\in\scrA$ and each
partition of $\cU$ into finitely many (equivalently, two) pieces, one of these pieces belongs to $\scrA$.
\end{defn}

\begin{lem}[\cite{QRT}]\label{Ramseyan}
Assume that $\scrA$ is Ramseyan. Then $\scrA\subseteq\Omega$.
\end{lem}
\begin{proof}
Assume that $\scrA$ is Ramseyan, and $\scrA\not\subseteq\Omega.$
Fix $\cU\in\scrA\setminus\Omega$, and a finite subset $F\subseteq
X$ such that $F$ is not contained in any $U\in\cU$. Since $\cU$ is
a cover of $X$, $|F|\geq 2$. For each $C\subsetneq F$,
let $\cU_C=\{U\in\cU : U\cap F=C\}$.
Then $\cU=\Union_{C\in P(F)\sm\{F\}}\cU_C$ is a partition of $\cU$ into finitely many pieces.

As $\scrA$ is Ramseyan, there is $C\subsetneq F$ such that
$\cU_C\in\scrA$. But then the elements of $F\sm C$ are not covered
by any member of $\cU_C$. A Contradiction.
\end{proof}

Examples of Ramseyan collections of covers are $\Omega$ and $\Gamma$, defined in the introduction.
We will give one more example in Section \ref{tau*}.

A cover $\cU$ of $X$ is \emph{multifinite} \cite{strongdiags} if there exists a
partition of $\cU$ into infinitely many finite covers of $X$.

\begin{defn}[The \emph{Gimel operator} on families of covers]
Let $\scrA$ be a family of covers of $X$. $\grpa$
is the family of all covers $\cU$ of $X$ such that:
Either $\cU$ is multifinite, or there exists a partition $\mathcal{P}$ of $\cU$
into finite sets such that
$\{\bigcup\mathcal{F}:\mathcal{F}\in\mathcal{P}\}\setminus\{X\}\in\scrA$.
\end{defn}

\begin{rem}
For each $\scrA$, $\scrA\sbst\grp{\scrA}$.
\end{rem}

An element of $\grpa$ will be called \emph{$\scrA$-glueable}. This
explains our choice of the Hebrew letter \emph{Gimel} ($\gimel$).

\begin{defn}
A cover $\cV$ is a \emph{finite-to-one derefinement} of a cover
$\cU$, if there exists a finite-to-one surjection $f:\cU\to \cV$
such that for each $U\in\cU$, $U\sbst f(U)$.

$\scrA$ is \emph{finite-to-one derefinable} if for each $\cU\in\scrA$ and
each finite-to-one derefinement $\cV\in\cO$ of $\cU$, $\cV\in\scrA$.
\end{defn}

A useful tool in the study of selection principles and their relation to
Ramsey theory is the game $\gfin(\scrA,\scrB)$. This game
is played by two players, ONE and TWO, with an inning per each
natural number $n$.
At the $n$th inning ONE chooses a cover $\cU_{n}\in\scrA$ and TWO chooses a finite
subset $\mathcal{F}_{n}$ of $\cU_{n}$. TWO wins if
$\bigcup_{n\in\mathbb{N}}\mathcal{F}_{n}\in\scrB$.
Otherwise, ONE wins.

Our goal in this section is proving the following.

\begin{thm}\label{main}
Let $\scrB$ be Ramseyan and finite-to-one derefinable.
The following are equivalent:
\begin{enumerate}
\item $\ufin(\cO,\scrB)$.
\item $\sfin(\cO,\cO)$ and $\Lambda=\grpb$.
\item ONE has no winning strategy in $\gfin(\Omega,\grpb)$.
\item $\sfin(\Omega,\grpb)$.
\end{enumerate}
Moreover, in (3) and (4), $\Omega$ can be replaced by any of
$\Lambda$ or $\Gamma$.
\end{thm}

We prove this theorem in a sequence of lemmas.
As it may be of independent interest,
some of these lemmas use weaker (or no) requirements on $\scrB$
than those posed in Theorem \ref{main}.

$\cU$ is a \emph{large cover} of $X$ if each point $x\in X$ belongs to
infinitely many $U\in\cU$. Let $\Lambda$ be the collection of all countable large covers of $X$.

The following proof is similar to that of \cite[Lemma 8]{coc7}.

\begin{lem}\label{large-is-B-glueable}
Assume that $\scrB$ is Ramseyan and finite-to-one derefinable.
Then $\ufin(\cO,\scrB)$ implies $\Lambda=\grpb$.
\end{lem}
\begin{proof}
By Lemma \ref{Ramseyan}, $\scrB\subseteq\Omega$. Thus, each
$\cU\in\grpb$ is large. Let $\cU$ be a large cover of $X$.
If $\cU$ is multifinite, then $\cU\in\grpb$, and we are done.

We now treat the remaining two cases.

\subsubsection*{Case 1.} $\cU$ has no finite subcover. Let $\{U_{n}:n\in\N\}$
bijectively enumerate a large cover $\cU$ of $X$. We may assume
that no finite subset of $\cU$ covers $X$: If $\cU$ contains
infinitely many disjoint finite subcovers then it is multifinite.

For $m,n\in\N$, we use the convenient notation
$$U_{[m,n)}=\Union_{m\le i < n}U_{i},$$
with the convention that $U_{[m,n)}=\emptyset$ whenever $n\le m$.
For each $n$, define $\cU_{n}=\{U_{[n,m)} : m\in\N\}$.
Each $\cU_{n}$ is a $\gamma$-cover of $X$, and in particular a cover of $X$.

Applying $\ufin(\cO,\scrB)$, choose for each $n$ a finite
$\cF_{n}\sbst \cU_{n}$ such that
$\{\bigcup\cF_{n}:n\in\N\}\in\scrB$.
For each $n$, there is $m_n\ge n$ such that $\Union\cF_n = U_{[n,m_n)}$.

Let $k_1 = 1$, and $k_2=m_1+1$.
Having defined $k_{n-1}$ and $k_n$, choose $k_{n+1}$ such that:
\be
\itm $k_{n+1}>m_1,m_2,\dots,m_{k_n}$;
\itm There is $i$ such that $k_n\le i<k_{n+1}$ and $U_{[i,m_i)}\neq\emptyset$; and
\itm $U_{[k_{n-1},k_{n+1})}\nin \{U_{[k_{i-1},k_{i+1})} : i<n\}$.
\ee
(3) is possible since $U_{[1,k_n)}\neq X$.

For each $n$, let $\cV_n =\{U_{[i,m_i)} : k_n\le i<k_{n+1}\}$.
As
$$\Union_n\cV_{2n-1}\cup\Union_n\cV_{2n}=\{U_{[i,m_i)} : i\in\N\}\in\scrB$$
and $\scrB$ is Ramseyan, there is $j\in\{0,1\}$ such that $\Union_n\cV_{2n-j}\in\scrB$.
We consider the case $j=0$ (the other case can be treated similarly).

For each $n$, each element of $\cV_{2n}$ has the form
$U_{[i,m_i)}$ with $k_{2n}\le i<k_{2n+1}$. By (1),
$U_{[i,m_i)}\sbst U_{[k_{2n},k_{2n+2})}$. Thus,
$\{U_{[k_{2n},k_{2n+2})} : n\in\N\}$ is a finite-to-one
derefinement of $\Union_n\cV_{2n}$, and is therefore a member of
$\scrB$. As $\scrB$ is finite-to-one derefinable,
$\{U_{[1,k_4)}\}\cup\{U_{[k_{2n},k_{2n+2})} : n>1\}\in\scrB$,
either, and this witnesses that the partition of $\cU$ into the
pieces $\{U_i : 1\le i < k_4\}$ and $\{U_i : k_{2n}\le i <
k_{2n+2}\}$, $n>1$, is as required in the statement $\cU\in\grpb$.

\subsubsection*{Case 2.} $\cU$ has only finitely many disjoint finite
subcovers. Let $\cF$ be the family of all elements in these finite subcovers.
$\cU\sm\cF$ is a large cover of $X$ not containing any finite subcover.
By what we have just proved, $\cU\sm\cF\in\grpb$.
As $\cU\sm\cF$ is not multifinite, there is a partition $\cP$ of $\cU\sm\cF$ into finite pieces,
such that $\{\Union\cV : \cV\in\cP\}\in\scrB$.
Fix $\cV_0\in\cP$. $\cP' = \{\cV_0\cup\cF\}\cup\cP\sm\{\cV_0\}$ is
a partition of $\cU$ into finite pieces.
Define $f:\{\Union\cV : \cV\in\cP\}\to \{\Union\cV : \cV\in\cP'\}$ by
$f(\Union\cV)=\Union(\cV\cup\cF)$ if $\Union\cV=\Union\cV_0$, and $f(\Union\cV)=\Union\cV$ otherwise.
As $f$ is finite-to-one and $\scrB$ is finite-to-one derefinable,
$\{\Union\cV : \cV\in\cP'\}\in\scrB$, and thus $\cU\in\grpb$.
\end{proof}

Note that for each family of covers $\scrB$,
$\ufin(\cO,\scrB)$ implies $\ufin(\cO,\cO)$. Clearly, $\ufin(\cO,\allowbreak\cO)=\sfin(\cO,\cO)$.
Thus, Lemma \ref{large-is-B-glueable} shows that the implication $(1)\Impl (2)$ of Theorem \ref{main}
holds for each Ramseyan and bijectively derefinable family of covers $\scrB$.

The following will be used often.

\begin{lem}[\cite{coc1}]\label{basic2}
$\sfin(\cO,\cO)=\sfin(\Lambda,\Lambda)=\sfin(\Omega,\Lambda)=\sfin(\Gamma,\Lambda)$.
\end{lem}

\begin{lem}\label{basicgames}
The following are equivalent:
\be
\itm $\sfin(\cO,\cO)$.
\itm ONE has no winning strategy in $\gfin(\Lambda,\Lambda)$.
\itm ONE has no winning strategy in $\gfin(\Omega,\Lambda)$.
\itm ONE has no winning strategy in $\gfin(\Gamma,\Lambda)$.
\ee
\end{lem}
\begin{proof}
Recall that $\sfin(\cO,\cO)=\sfin(\Lambda,\Lambda)$.
$\sfin(\Lambda,\Lambda)$ is equivalent to (2) \cite[Theorem 5]{OpPar}.

As $\Gamma\sbst\Omega\sbst\Lambda$, $(2)\Impl (3)\Impl (4)$.
But $(4)$ implies $\sfin(\Gamma,\Lambda)$, which is the same as (1) \cite{coc1}.
\end{proof}

\begin{cor}
The conjunction of $\sfin(\cO,\cO)$ and $\Lambda=\grpb$ implies
that ONE has no winning strategy in any of the games
$\gfin(\Lambda,\grpb)$, $\gfin(\Omega,\grpb)$, or $\gfin(\Gamma,\grpb)$.
\end{cor}
\begin{proof}
Lemma \ref{basicgames} and the assumption $\Lambda=\grpb$.
\end{proof}

This gives $(2)\Impl (3)$ of Theorem \ref{main}.
$(3)\Impl (4)$ in that theorem is clear.
It remains to show that $(4)\Impl (1)$.
As $\Gamma\sbst\Omega\sbst\Lambda$, it suffices to prove the following.

\begin{lem}\label{lem:equiv-sfin-ufin2}
Assume that $\scrB$ is finite-to-one derefinable.
Then $\sfin(\Gamma,\grpb)$ implies $\ufin(\cO,\scrB)$.
\end{lem}
\begin{proof}
Assume that $X$ satisfies $\sfin(\Gamma,\grpb)$.
As $\ufin(\cO,\scrB)=\ufin(\Gamma,\scrB)$ \cite{coc1},
it suffices to prove that $X$ satisfies $\ufin(\Gamma,\scrB)$.

Let $\cU_{n}$, $n\in\N$, be disjoint open $\gamma$-covers of $X$
which do not contain a finite subcover. Enumerate each $\cU_{n}$
bijectively as $\{U_{k}^{n}:k\in\N\}$. For each $n$, let
$$\cV_{n}=\{U_{m}^{1}\cap U_{m}^{2}\cap\dots\cap U_{m}^{n}:m\in\N\}.$$
For each $n$, $\cV_{n}$ is an open $\gamma$-cover of $X$. Apply
$\sfin(\Gamma,\grpb)$ to obtain for each $n$ a finite subset
$\mathcal{F}_{n}\subseteq\cV_{n}$ such that
$\bigcup_{n}\mathcal{F}_{n}\in\grpb$. Each
$U\in\bigcup_{n}\mathcal{F}_{n}$ is a subset of some element of
$\cU_{1}$, hence for each finite subset
$\mathcal{F}\subseteq\bigcup_{n}\mathcal{F}_{n}$,
$\bigcup\mathcal{F}\neq X$. Therefore $\bigcup_{n}\mathcal{F}_{n}$
is not multifinite.

Let $\{\mathcal{X}_{m}:m\in\N\}$ be a partition of
$\bigcup_{n}\mathcal{F}_{n}$ into finite pieces such that
$\{\bigcup\mathcal{X}_{m}:m\in\N\}\in\scrB$. Let\[
f(m)=\min\{ k:\mathcal{X}_{m}\cap\mathcal{F}_{k}\neq\emptyset\},\]
and put \[
\mathcal{Y}_{n}=\bigcup_{m\in f^{-1}(n)}\mathcal{X}_{m}\] The sets
$\{ f^{-1}(n):n\in\N\}$ form a partition of $\N$.
Since each $\mathcal{F}_{k}$ is finite and the $\mathcal{X}_{m}$'s
are disjoint, $f^{-1}(n)$ is finite for all $n$. It follows that
each $\mathcal{Y}_{n}$ is a finite set. Each member of
$\mathcal{Y}_{n}$ belong to some $\mathcal{F}_{k}\subseteq\cV_{k}$
for some $k\geq n$.

For each $n$, choose $\psi(n)\in\N$ such that:
\begin{enumerate}
\item If $U_{k}^{n}\in\cU_{n}$ appear as term in the sets of
$\mathcal{Y}_{n}$ then $\psi(n)\geq k$. \item The sets
$\bigcup_{k\leq\psi(n)}U_{k}^{n}$ are distinct for different
values of $n$.
\end{enumerate}
This is possible since $\{\bigcup_{k\leq
m}U_{k}^{n}:m\in\N\}$, $n\in\N$, are $\gamma$-covers.

Define $\mathcal{Z}_{n}=\{ U_{k}^{n}:k\leq\psi(n)\}$. The sets $\mathcal{Z}_{n}$
are finite and disjoint. For each $n\in\N$, $\bigcup\mathcal{Y}_{n}\subseteq\bigcup\mathcal{Z}_{n}\neq X$.
Hence $\{\bigcup\mathcal{Z}_{n}:n\in\N\}$ is a finite derefinement
of $\{\bigcup\mathcal{X}_{n}:n\in\N\}$. Therefore $\{\bigcup\mathcal{Z}_{n}:n\in\N\}\in\scrB$
and the sequence $\{\mathcal{Z}_{n}\}_{n\in\N}$ witnesses
that $X$ has the property $\ufin(\Gamma,\scrB)$.
\end{proof}

This completes the proof of Theorem \ref{main}.

\section{Partition relations for glueable covers}

The symbol $[A]^{n}$ denotes the set of $n$-element subsets of $A$.
For a positive integer $k$, the \emph{Baumgartner-Taylor partition relation} \cite{BT}
\[
\scrA\to\lceil\scrB\rceil^{2}_{k}
\]
denotes the following statement:
For each $A$ in $\scrA$ and each $f:[A]^{2}\to\{1,\dots,k\}$,
there are
\be
\itm $B\subseteq A$ such that $B\in\scrB$;
\itm A partition of $B$ into finite pieces $B=\bigcup_{n\in\N} B_n$; and
\itm $j\in\{1,\dots,k\}$,
\ee
such that $f(\{U,V\})=j$ for all $U,V\in B$ which do not belong to the same $B_n$.

The Baumgartner-Taylor partition relation is one of the most important partition
relations in the studies of open covers and their combinatorial properties --
see \cite{KocRamsey} for a survey of this field.

\begin{lem}[\cite{QRT}]\label{lem:B-infinite-then-A-in-Omega}
If each member of $\scrB$ is infinite and $\scrA\to\lceil\scrB\rceil_{2}^{2}$ holds,
then $\scrA\subseteq\Omega$.
\end{lem}

Together with Theorem \ref{main}, the following gives a Ramsey theoretic characterization
of properties of the form $\ufin(\cO,\scrB)$.

\begin{thm}\label{main2}
Assume that $\scrB$ is Ramseyan and finite-to-one derefinable.
The following are equivalent:
\be
\itm $\sfin(\Omega,\grpb)$.
\itm For each $k$, $\Omega\to\lceil\grpb\rceil_{k}^{2}$ holds.
\itm $\Omega\to\lceil\grpb\rceil_{2}^{2}$.
\ee
\end{thm}
\begin{proof}
$(1\Impl 2)$ This follows from Theorem \ref{main}
and the following.
\begin{lem}[\cite{coc7}]
Assume that $\scrA$ is Ramseyan. If ONE has no winning strategy
in the game $\gfin(\scrA,\scrB)$, then for each $k$, $\scrA\to\lceil\scrB\rceil^{2}_{k}$ holds.
\end{lem}

$(3\Impl 1)$ Assume that $X$ satisfies $\Omega\to\lceil\grpb\rceil_{2}^{2}$.
By Theorem \ref{main}, it suffices to show that $X$ satisfies $\ufin(\Gamma,\scrB)$.

Let
$\cU_{n}$, $n\in\N$, be open $\gamma$-covers of $X$ which
do not contain a finite subcover. Enumerate each $\cU_{n}$
bijectively as $\{U_{k}^{n}:k\in\N\}$. For each $n$,
define
\[
\cV_{n}=\{ U_{k}^{1}\cap U_{k}^{2}\cap\dots\cap U_{k}^{n}:k\in\N\},
\]
and let
$\cV=\bigcup_{n\in\N}\cV_{n}$. Then, $\cV$ is an $\omega$-cover of $X$.
For each element of $\cV$ fix a representation of
the form $U_{k}^{1}\cap U_{k}^{2}\cap\dots\cap U_{k}^{n}$.
Define a function $f:[\cV]^{2}\to\{1,2\}$ by\[ f(\{
V_{1},V_{2}\})=\begin{cases}
1 & \mbox{if }V_{1}\mbox{ and }V_{2}\mbox{ are from the same }\cV_{n},\\
2 & \mbox{otherwise.}\end{cases}\]

Choose $\cW\subseteq\cV$ such that $\cW\in\grpb$,
a partition $\cW=\Union_k\cW_k$ into finite pieces, and a color
$j\in\{1,2\}$, such that for $A$ and $B$ from distinct
$\cW_{k}$'s, $f(\{ A,B\})=j$.
Consider the possible values of $j$.

$j=1$: Then there is an $n$ such that for all
$A\in\cW$ we have $A\subseteq U_{n}^{1}\neq X$. Hence
$\cW$ is not a cover. Contradiction.

$j=2$: Let $\mathcal{F}_{n}=\cW\cap\cV_{n}$. Then each
$\mathcal{F}_{n}$ is finite. From this point, the proof continues
as in the proof of Lemma \ref{lem:equiv-sfin-ufin2}.
\end{proof}

\section{Selecting one element from each cover}

We now consider the following selection principle.
\begin{description}
\item [{$\sone(\scrA,\scrB)$}] For each sequence
$\{\cU_{n}\}_{n\in\N}$ of elements of $\scrA$, there
exist $U_{n}\in\cU_{n}$, $n\in\N$, such that $\{U_{n}:n\in\N\}\in\scrB$.
\end{description}

The corresponding game $\gone(\scrA,\scrB)$, is defined as follows:
At the $n$th inning ONE chooses a cover $\cU_{n}\in\scrA$ and TWO chooses
$U_{n}\in\cU_{n}$. TWO wins if
$\{U_{n}:n\in\N\}\in\scrB$. Otherwise, ONE wins.

The corresponding partition relation, called the ordinary partition relation, is defined as follows.
For positive integers $n$ and $k$,
\[
\scrA\to(\scrB)_{k}^{n}
\]
means:
For each $A\in\scrA$ and each $f:[A]^{n}\to\{1,\dots,k\}$,
there is $B\subseteq A$ such that $B\in\scrB$, and $f|_{[B]^{n}}$ is constant.

The following theorem was proved in \cite{coc7} for $\scrB=\Gamma$,
and in \cite{coc8} for $\scrB=\Omega$.

\begin{thm}\label{main3}
Let $\scrB$ be Ramseyan and finite-to-one derefinable.
The following are equivalent.
\begin{enumerate}
\item $\sone(\cO,\cO)$ and $\ufin(\cO,\scrB)$.
\item $\sone(\Lambda,\grpb)$.
\item $\sone(\Omega,\grpb)$.
\item ONE has no winning strategy in the game $\gone(\Omega,\grpb)$.
\item $\Omega\to(\grpb)_{2}^{2}$.
\item $\Omega\to(\grpb)_{k}^{2}$ for all $k$.
\end{enumerate}
\end{thm}
\begin{proof}
$(1\Impl 2)$
$\sone(\cO,\cO)=\sone(\Lambda,\Lambda)$.
By Lemma \ref{large-is-B-glueable}, $\Lambda=\grpb$ for $X$. Thus, $X$ satisfies $\sone(\Lambda,\grpb)$.

$(2\Impl 3)$ $\Omega\sbst\Lambda$.

$(3\Impl 1)$ As $\sone(\Omega,\grpb)$ implies $\sfin(\Omega,\allowbreak\grpb)$,
we have by Lemma \ref{lem:equiv-sfin-ufin2} that $\ufin(\cO,\scrB)$ holds,
and that $\grpb=\Lambda$.
Thus, $X$ satisfies $\sone(\Omega,\Lambda)$, which is the same as $\sone(\cO,\cO)$ \cite{coc1}.

$(1\Impl 4)$ By \cite[Theorem 3]{OpPar}, $\sone(\cO,\cO)$
implies that ONE does not have a strategy in $\gone(\Lambda,\Lambda)$,
and in particular in $\gone(\Omega,\Lambda)$.
Again, use Theorem \ref{main} to get that $\Lambda=\grpb$.

$(4\Impl 6)$ Follows from \cite[Theorem 1]{coc7}.

$(6\Impl 5)$ is immediate.

$(5\Impl 3)$ As $\scrB$ is Ramseyan, $\scrB\sbst\Omega$, and therefore $\grp{\scrB}\sbst\Lambda$.
Thus, (5) implies $\Omega\to(\Lambda)_{k}^{2}$. Using the methods of \cite{coc7}, one can
prove that $\Omega\to(\Lambda)_{k}^{2}$ implies $\sone(\Omega,\Lambda)$ \cite{QRT}.
Clearly, (5) also implies $\Omega\to\lceil\grpb\rceil_{2}^{2}$, and by Theorem \ref{main2},
we get $\Lambda=\grpb$.
\end{proof}

\section{Applications}\label{tau*}

\subsection{$\gamma$-covers}

As every infinite subset of a $\gamma$-cover is again a $\gamma$-cover of the same space, $\Gamma$ is Ramseyan.

\begin{lem}
$\Gamma$ is finite-to-one derefinable.
\end{lem}
\begin{proof}
Assume that $\cU\in\Gamma$ and $f:\cU\to\cV$ is
finite-to-one and surjective. As $f$ is finite-to-one and $\cU$ is infinite, $\cV$ is infinite.
Assume that $x\in X$ and $\cW = \{V\in \cV : x\nin V\}$ is infinite.
For each $V\in\cW$ and each $U\in f^{-1}(V)$, $U\sbst V$ and thus $x\nin U$.
As $f$ is surjective, $\Union_{V\in\cW} f^{-1}(V)$ is infinite. A contradiction.
\end{proof}

Thus, we can directly apply Theorems \ref{main}, \ref{main2}, and \ref{main3},
and obtain the following.

\begin{thm}[\cite{coc7}]
The following are equivalent:
\begin{enumerate}
\item $\ufin(\cO,\Gamma)$.
\item $\sfin(\cO,\cO)$ and $\Lambda=\grp{\Gamma}$.
\item ONE has no winning strategy in $\gfin(\Omega,\grp{\Gamma})$.
\item $\sfin(\Omega,\grp{\Gamma})$.
\itm For each $k$, $\Omega\to\lceil\grp{\Gamma}\rceil_{k}^{2}$ holds.
\itm $\Omega\to\lceil\grp{\Gamma}\rceil_{2}^{2}$.\hfill\qed
\end{enumerate}
\end{thm}

\begin{thm}[\cite{coc7}]
The following are equivalent.
\begin{enumerate}
\item $\sone(\cO,\cO)$ and $\ufin(\cO,\Gamma)$.
\item $\sone(\Lambda,\grp{\Gamma})$.
\item $\sone(\Omega,\grp{\Gamma})$.
\item ONE has no winning strategy in the game $\gone(\Omega,\grp{\Gamma})$.
\item $\Omega\to(\grp{\Gamma})_{2}^{2}$.
\item $\Omega\to(\grp{\Gamma})_{k}^{2}$ for all $k$.\hfill\qed
\end{enumerate}
\end{thm}

\subsection{$\omega$-covers}

\begin{defn}
A cover $\cV$ is a \emph{derefinement} of a cover
$\cU$ if $\cU$ refines $\cV$.
$\scrA$ is \emph{derefinable} if for each $\cU\in\scrA$ and
each derefinement $\cV\in\cO$ of $\cU$, $\cV\in\scrA$.
\end{defn}

$\Omega$ is derefinable, and in particular finite-to-one derefinable.

\begin{lem}[folklore]\label{folk}
$\Omega$ is Ramseyan.
\end{lem}
\begin{proof}
Assume that $\cU\in\Omega$ and
$\cU = \cU_1\cup\dots\cU_n$ and no $\cU_i\in\Omega$. For each $i$, choose
a finite subset $F_i$ of $X$ witnessing $\cU_i\nin\Omega$. Then $F=F_1\cup\dots\cup F_n$
is not covered by any element of $\cU$. A contradiction.
\end{proof}

In the forthcoming Theorem \ref{coc7a}, we reproduce the
statements of Theorems 2 and 3 of \cite{coc8}. One direction in
the proof of Theorem 3 in \cite{coc8} uses Theorem 4 of
\cite{coc7}, which in turn requires that $\grp{\Omega}$ is
derefinable. Unfortunately, by Theorem 2 of \cite{coc8}, spaces
dealt with in this theorem only have $\Lambda=\grp{\Omega}$. But
$\Lambda$ is not derefinable: Fix distinct $a,b,x_n\in X$,
$n\in\N$. Then the large cover $\{X\sm\{a,x_n\},X\sm\{b,x_n\} :
n\in\N\}$ refines $\{X\sm\{a\},X\sm\{b\}\}$. Our results give a
corrected proof of this direction.

\begin{thm}\label{coc7a}
The following are equivalent:
\begin{enumerate}
\item $\ufin(\cO,\Omega)$.
\item $\sfin(\cO,\cO)$ and $\Lambda=\grp{\Omega}$.
\item ONE has no winning strategy in $\gfin(\Omega,\grp{\Omega})$.
\item $\sfin(\Omega,\grp{\Omega})$.
\itm For each $k$, $\Omega\to\lceil\grp{\Omega}\rceil_{k}^{2}$ holds.
\itm $\Omega\to\lceil\grp{\Omega}\rceil_{2}^{2}$.
\end{enumerate}
\end{thm}
\begin{proof}
$\Omega$ is derefinable and Ramseyan (Lemma \ref{folk}).
Apply Theorems \ref{main} and \ref{main2}.
\end{proof}

As in the previous theorem, the following Theorem \ref{coc7b}
reproduces Theorem 5 of \cite{coc8} and fixes a problem similar
to the above-mentioned one in the original proof of the implication $(5)\Impl (3)$ below.

\begin{thm}\label{coc7b}
The following are equivalent.
\begin{enumerate}
\item $\sone(\cO,\cO)$ and $\ufin(\cO,\Omega)$.
\item $\sone(\Lambda,\grp{\Omega})$.
\item $\sone(\Omega,\grp{\Omega})$.
\item ONE has no winning strategy in the game $\gone(\Omega,\grp{\Omega})$.
\item $\Omega\to(\grp{\Omega})_{2}^{2}$.
\item $\Omega\to(\grp{\Omega})_{k}^{2}$ for all $k$.
\end{enumerate}
\end{thm}
\begin{proof}
Apply Theorem \ref{main3}.
\end{proof}

\subsection{$\tau^*$-covers}
Let $\roth=\{A\sbst\N : |A|=\aleph_0\}$. For $A,B\in\roth$, $A\as B$ means
that $A\sm B$ is finite.
A family $Y\subseteq[\N]^{\aleph_{0}}$ is
\emph{linearly refinable} if for each $y\in Y$ there exists an infinite subset
$\hat{y}\subseteq y$ such that the family $\hat{Y}=\{\hat{y}:y\in Y\}$
is linearly ordered by $\subseteq^{*}$.

A countable cover $\cU=\{U_{n}:n\in\N\}$ of $X$ is a
\emph{$\tau^{*}$-cover} of $X$ if $\left\{ \left\{ n:x\in U_{n}\right\} :x\in X\right\}$
is linearly refinable. $\Tau^*$ is the collection of all $\tau^{*}$-covers.
$\Gamma\sbst\Tau^*\sbst\Omega$.

$\Tau^*$ is derefinable \cite{tautau}.
In particular, $\Tau^*$ is finite-to-one derefinable.
(This latter assertion is easier to see.)

\begin{prop}\label{lem:T-star-persistent}
$\Tau^*$ is Ramseyan.
\end{prop}
\begin{proof}
Let $\{ U_{n}:n\in\N\}$ be a bijective enumeration of a
$\tau^{*}$-cover $\cU$ of $X$.
For each $x\in X$, let $x_{\cU}=\{ n:x\in U_{n}\}$, and let $\hat{x}_{\cU}$ be an
infinite subset of $x_{\cU}$ such that the sets $\hat{x}_{\cU}$ are linearly ordered by
$\subseteq^{*}$.

Consider a partition $\cU=\cV\cup(\cU\setminus\cV)$.
Define $A=\{ n:U_{n}\in\cV\}$. We may assume that both $A$ and its
complement are infinite.

For each $x\in X$, define $\hat{x}_{\cV}=\hat{x}_{\cU}\cap A$ and
$\hat{x}_{\cU\setminus\cV}=\hat{x}_{\cU}\cap A^{c}$. If
$\{\hat{x}_{\cV}:x\in X\}\subseteq\roth$ or
$\{\hat{x}_{\cU\setminus\cV}:x\in
X\}\subseteq\roth$ then we are done. If this
is not the case, then there are some $x,y\in X$ such that
$\hat{x}_{\cV}$ and $\hat{y}_{\cU\setminus\cV}$ are finite.
Without loss of generality, assume that
$\hat{y}_{\cU}\subseteq^{*}\hat{x}_{\cU}$. Thus,\[
\hat{y}_{\cV}=\hat{y}_{\cU}\cap A\subseteq^{*}\hat{x}_{\cU}\cap
A=\hat{x}_{\cV}\] but $\hat{y}_{\cV}$ is infinite and
$\hat{x}_{\cV}$ is finite. A contradiction.
\end{proof}

By Theorems \ref{main} and \ref{main2}, we have the following.

\begin{thm}
The following are equivalent:
\begin{enumerate}
\item $\ufin(\cO,\Tau^*)$.
\item $\sfin(\cO,\cO)$ and $\Lambda=\grp{\Tau^*}$.
\item ONE has no winning strategy in $\gfin(\Omega,\grp{\Tau^*})$.
\item $\sfin(\Omega,\grp{\Tau^*})$.
\itm For each $k$, $\Omega\to\lceil\grp{\Tau^*}\rceil_{k}^{2}$ holds.
\itm $\Omega\to\lceil\grp{\Tau^*}\rceil_{2}^{2}$.\hfill\qed
\end{enumerate}
\end{thm}

By Theorem \ref{main3}, we have the following.

\begin{thm}
The following are equivalent.
\begin{enumerate}
\item $\sone(\cO,\cO)$ and $\ufin(\cO,\Tau^*)$.
\item $\sone(\Lambda,\grp{\Tau^*})$.
\item $\sone(\Omega,\grp{\Tau^*})$.
\item ONE has no winning strategy in the game $\gone(\Omega,\grp{\Tau^*})$.
\item $\Omega\to(\grp{\Tau^*})_{2}^{2}$.
\item $\Omega\to(\grp{\Tau^*})_{k}^{2}$ for all $k$.\hfill\qed
\end{enumerate}
\end{thm}

\subsection*{Acknowledgment} We thank the referee for the useful comments on this paper.


\begin{thebibliography}{99}

\bibitem{coc8}
L.\ Babinkostova, Lj.\ Ko\v{c}inac, and M.\ Scheepers,
\emph{Combinatorics of open covers (VIII)},
Topology and its Applications \textbf{140} (2004), 15--32.

\bibitem{BT}
J. Baumgartner and A. Taylor,
\emph{Partition theorems and ultrafilters}, \textbf{Transactions of the American Mathematical
Society} 241 (1978), 283-309.

\bibitem{Hure25}
W. Hurewicz,
\emph{\"Uber eine Verallgemeinerung des Borelschen Theorems},
Mathematische Zeitschrift \textbf{24} (1925), 401--421.

\bibitem{Hure27}
W. Hurewicz,
\emph{\"Uber Folgen stetiger Funktionen},
Fundamenta Mathematicae \textbf{9} (1927), 193--204.

\bibitem{KocRamsey}
Lj. Ko\v{c}inac,
\emph{Generalized Ramsey theory and topological properties: A survey},
\textbf{Rendiconti del Seminario Matematico di Messina}, Serie II, 9 (2003), 119--132.

\bibitem{coc7}
Lj. Ko\v{c}inac and M. Scheepers,
\emph{Combinatorics of open covers (VII): Groupability},
Fundamenta Mathematicae \textbf{179} (2003), 131--155.

\bibitem{Menger24}
K. Menger,
\emph{Einige \"Uberdeckungss\"atze der Punktmengenlehre},
Sitzungsberichte der Wiener Akademie \textbf{133} (1924), 421--444.

\bibitem{QRT}
N. Samet and B. Tsaban,
\emph{Ramsey theory of open covers},
in progress.

\bibitem{coc1}
M. Scheepers,
\emph{Combinatorics of open covers I: Ramsey theory},
Topology and its Applications \textbf{69} (1996), 31--62.

\bibitem{OpPar}
M. Scheepers,
\emph{Open covers and partition relations},
Proceedings of the American Mathematical Society \textbf{127} (1999), 577--581.

\bibitem{tautau}
B. Tsaban,
\emph{Selection principles and the minimal tower problem},
Note di Matematica \textbf{22} (2003), 53--81.

\bibitem{strongdiags}
B. Tsaban,
\emph{Strong $\gamma$-sets and other singular spaces},
Topology and its Applications \textbf{153} (2005), 620--639.

\end{thebibliography}
\end{document}